\documentclass[a4paper,reqno]{amsart}

\usepackage[utf8]{inputenc}
\usepackage[T1]{fontenc}
\usepackage{times,eulervm}
\usepackage{subcaption}
\usepackage{enumerate}
\usepackage{booktabs}
\usepackage{hyphenat}
\usepackage{amssymb}
\usepackage{siunitx}
\usepackage{bm}
\usepackage{hyperref}
\usepackage[capitalize]{cleveref}

\usepackage{ao-abbr}
\usepackage[final]{ao-ximg}
\usepackage{ao-matrix}
\usepackage{ao-math-opt}
\usepackage{ao-math-fields}
\usepackage{ao-math-symbols}
\usepackage{ao-code-names}

\newcommand\ff{\bm f}
\newcommand\spext[1][]{\tsp{ext\ifgiven{#1}{,#1}{}}}
\newcommand\spint{\tsp{int}}
\newcommand\bh{\bm h}
\newcommand\blam{\bm\lambda}
\newcommand\bl{\bm l}
\newcommand\bq{\bm q}
\newcommand\br{\bm r}
\newcommand\bs{\bm s}
\newcommand\bv{\bm v}

\newcommand\bx{\bm x}
\newcommand\by{\bm y}
\newcommand\bB{\bm B}
\newcommand\bF{\bm F}
\newcommand\bH{\bm H}
\newcommand\bK{\bm K}
\newcommand\bOmega{\bm\Omega}

\newcommand\incr[1][k]{\bm p_{#1}}
\newcommand\incN[1][k]{\bm p_{#1}^N}
\newcommand\faer{\ff\tsp{aero}}
\newcommand\Maer{\bK\tsp{aero}}
\newcommand\Mstr{\bK\tsp{str}}

\newcommand\IfAM{Leibniz Universität Hannover\\
  Institute of Applied Mathematics\\
  Wel\-fen\-gar\-ten 1\\30167 Hannover\\Germany}
\newcommand\ISD{Leibniz Universität Hannover\\
  Institute of Statics and Dynamics\\
  Appelstr.\ 9a\\30167 Hannover\\Germany}
\newcommand\orcid[1]{\hfill\break\hbox{}\kern\parindent ORCID #1}

\begin{document}

\title
[Accelerating Aeroelastic UVLM Simulations by Inexact Newton Algorithms]
{Accelerating Aeroelastic UVLM Simulations\\by Inexact Newton Algorithms}

\author[J. Schubert]{Jenny Schubert}
\address{Jenny Schubert\\\IfAM\orcid{0000-0001-5857-4689}}
\email{schubert@ifam.uni-hannover.de}
\urladdr{ifam.uni-hannover.de/schubert}

\author[M. C. Steinbach]{Marc C. Steinbach}
\address{Marc C. Steinbach\\\IfAM\orcid{0000-0002-6343-9809}}
\email{mcs@ifam.uni-hannover.de}
\urladdr{ifam.uni-hannover.de/mcs}

\author[C. Hente]{Christian Hente}
\address{Christian Hente\\\ISD\orcid{0009-0003-2229-6770}}
\email{c.hente@isd.uni-hannover.de}

\author[D. Märtins]{David Märtins}
\address{David Märtins\\\ISD\orcid{0009-0009-2395-4126}}
\email{d.maertins@isd.uni-hannover.de}
\urladdr{isd.uni-hannover.de/en/institute/team/david-maertins-msc}

\author[D. Schuster]{Daniel Schuster}
\address{Daniel Schuster\\\ISD\orcid{0000-0002-0072-9614}}
\email{d.schuster@isd.uni-hannover.de}
\urladdr{isd.uni-hannover.de/en/schuster}

\begin{abstract}
  We consider the aeroelastic simulation of flexible mechanical structures
  submerged in subsonic fluid flows at low mach numbers.
  The nonlinear kinematics of flexible bodies
  are described in the total Lagrangian formulation
  and discretized by finite elements.
  The aerodynamic loads are computed using the unsteady vortex-lattice method
  wherein a free wake is tracked over time.
  Each implicit time step in the dynamic simulation then requires
  solving a nonlinear equation system in the structural variables
  with additional aerodynamic load terms.
  Our focus here is on the efficient numerical solution of this system
  by accelerating the Newton algorithm.
  The particular struture of the aeroelastic nonlinear system
  suggests the structural derivative as an approximation
  to the full derivative in the linear Newton system.
  We investigate and compare two promising algorithms
  based in this approximation, a quasi-Newton type algorithm
  and a novel inexact Newton algorithm.
  Numerical experiments are performed on a flexible plate and on a wind turbine.
  Our computational results show that the approximations
  can indeed accelerate the Newton algorithm substantially.
  Surprisingly, the theoretically preferable inexact Newton algorithm
  is much slower than the quasi-Newton algorithm,
  which motivates further research to speed up derivative evaluations.
\end{abstract}

\keywords{%
  Aeroelastic simulation,
  unsteady vortex-lattice method,
  implicit time integration,
  inexact Newton algorithms.
}

\subjclass[2010]{%
  49M15, 
  90C53, 
  74F10, 
  76B47
}

\date\today

\maketitle

\section{Introduction}
In a variety of applications, the subsonic flow around lifting surfaces separates at a defined location and is otherwise attached to the surface. These lifting surfaces are often formed by flexible bodies, for which the aerodynamic forces exerted are substantial enough to induce motional interactions with the structure. In these cases, the aeroelastic investigation of both structural and flow behavior is of particular interest.

A diverse range of applications highlights the importance of investigating the aeroelastic behavior within these specific fluid flow contexts, for instance helicopter rotors \cite{Wie_2009}, morphing wings \cite{Fonzi_2020}, and notably, the rotor blades of wind turbines. Recent trends show an increase in the size of wind turbine rotor blades, leading to higher flexibility \cite{OLAF_2020}. The high complexity of these systems precludes the use of closed analytical methods. Instead, a partitioned methodology is employed, involving an unsteady three-dimensional flow solver in conjunction with a nonlinear structural solver to ensure a comprehensive analysis of the structural system \cite{Wang_2016}. The choice of approach varies according to specific requirements. In the design phase, where numerous calculations are necessary, low-fidelity methods such as Blade Element Momentum Theory (BEMT) are often utilized, particularly in wind energy systems \cite{Ramos-Garcia_2022}. Conversely, for detailed investigations of fluid-structure interactions in a limited set of simulations, high-fidelity methods are preferred. These methods, which solve the three-dimensional Navier-Stokes equations using techniques like the Finite Volume Method, offer superior accuracy. However, they are associated with substantial computational demands, especially in scenarios involving variable fluid-structure boundaries \cite{Verstraete_2023}.

The structure itself is usually modelled numerically using finite elements, \eg, finite beam elements in cases where the overall structural dynamics are of interest, or finite shell models in cases where a detailed investigation of the local structural behavior such as the stress field is carried out \cite{Wang_2016}.

In many cases, mid-fidelity methods based on potential flow theory offer an attractive compromise between accurately predicting the unsteady three-dimensional flow behavior and feasible computational costs. One such method, which has been applied to investigate aeroelastic phenomena, is the unsteady vortex-lattice method (UVLM). Among others, Palacios et al.\ developed and verified a coupled UVLM-beam solver for the investigation of flexible aircraft \cite{Palacios:2010}. Murua and colleagues similarly investigated flight dynamics and aeroelastic effects on aircraft using a coupled UVLM-beam approach with linearization in state-space \cite{Murua:2012}. Hang et al.\ also performed an analytical sensitivity analysis of UVLM to investigate aircraft but coupled this with both linear beams and shells \cite{Hang:2020}. An analytical sensitivity analysis of UVLM was also performed by Stanford and Beran to optimize flapping wings \cite{Stanford:2010}. Mu\~noz-Simón et al.\ used UVLM to investigate the aerodynamic behavior of wind turbines and compared it with BEMT results, finding that UVLM better captures the unsteady effects occurring in real operating conditions \cite{Munoz_Simon:2020}. Ritter and colleagues validate the UVLM, including linearization in state-space, against experiments, using a highly flexible wing called Pazy Wing \cite{Ritter:2021}. Goizueta et al.\ compare different models of this wing to predict flutter as preparation for wind tunnel tests \cite{Goizueta:2021}. The Pazy Wing benchmark was used for validation in a comparison of several nonlinear aeroelastic solvers, especially with regard to flutter characteristics \cite{Ritter:2024}. The main difference of the aeroelastic framework \cite{Hente_et_al:2023} used within our work compared to the ones describe above is that we directly solve the nonlinear aeroelastic equlibrium without using a linearized state-space. The linearization of the aerodynamic forces is based on the unsteady Bernoulli equation with respect to generalized structural coordinates and velocities. We apply an implicit time integration scheme based on discrete time derivatives. For details on the aeroelastic framework, its UVLM implementation and linearization we refer the reader to \cite{Hente_et_al:2023}.
While Newton's method is the standard nonlinear system solver
in any implicit time integration,
inexact Newton algorithms depend on the application-specific structure.
To the best of our knowledge, our work is the first investigation
of inexact Newton algorithms for aeroelastic simulations using the UVLM.

The coupled problem can either be solved in a weak form, where only the influence of the aerodynamic forces on the structure is considered at a certain time, or in a strong form, where both problems are evaluated at the same time instance and the coupled nonlinear governing equation of the aeroelastic system has to be solved \cite{Lamei_2020}.
For systems experiencing nonlinear behavior, \eg, due to large deformations, a strong coupling is preferable. To integrate the governing equation in time, implicit time integration has proven more robust \cite{Hente_et_al:2023}. Consequently, the coupled nonlinear governing equation needs to be solved in each time-step. This is typically done by applying the Newton algorithm, which requires the linearization of the equation, \ie, the evaluation of a derivative matrix of a nonlinear function.

The main drawback of the (exact) Newton algorithm
is the potentially high cost
for evaluating the derivative matrix and
for solving the linear Newton system
at each iteration.
On the other hand, due to the excellent local convergence,
it is often possible to save computational cost
by solving the exact Newton system approximately
or by solving a system with an approximated derivative marix.
These two approaches are both promising in our context
because the full derivative matrix is the sum
of a structural derivative (which is cheap to compute and to invert)
and an aerodynamic derivative
(which has entries of smaller magnitude but is expensive to compute).
Then the first possiblity is a quasi-Newton type algorithm
wherein the approximate Newton step is computed
by neglecting the expensive aerodynamic derivative entirely,
\ie, solving the linear system with the structural inverse.
This algorithm has already been proposed in \cite{Hente_et_al:2023}.
The second possiblity is a so-called inexact Newton algorithm
wherein the same structural inverse is used
to solve the exact Newton system to a pre-specified accuracy
by iterative refinement.
This algorithm is novel for aeroelastic UVLM simulation.
Both approaches have their pros and cons.
In particular, the quasi-Newton type algorithm
is faster per iteration but may converge slowly or even diverge
whereas the inexact Newton algorithm is slower per iteration
but more robust and can even guarantee quadratic local convergence;
for details see \cref{sec:solution}.
The precise computational cost of each algorithm
depends on the individual costs of evaluating each of the two derivatives
and of solving with the full inverse versus the approximate inverse.
As it turns out, the overall situation is quite peculiar in our context.

The remainder of this paper is structured as follows.
In \cref{sec:simulation} we give a brief overview of the structural
and aerodynamic simulation models along with their coupling
by the unsteady vortex-lattice method.
Exact and approximate Newton algorithms and their specific application
within our aeroelastic simulation are presented in \cref{sec:solution}.
In \cref{sec:results} we conduct numerical experiments
with the three Newton type algorithms
on a flexible plate and on a wind turbine.
Finally we summarize our findings in \cref{sec:conclusion}.


\section{Aeroelastic Simulation}
\label{sec:simulation}
\subsection{Continuous governing equation}
This subsection provides a brief overview of the structural model employed in this article. For a more in-depth understanding and comprehensive explanations, we direct the reader to refer to the articles \cite{gebhardt_nonlinear_2017,  gebhardt_new_2020, Gebhardt_et_al:2020b, hente_modal_2019, hente_static_2021}.

The structural model relies on a rotation-free multibody system formalism within the Finite Element Method (FEM) presented in \cite{Gebhardt_et_al:2020b}, consisting of rigid bodies, geometrically exact beams and geometrically exact solid-degenerate shells \cite{gebhardt_nonlinear_2017}. The kinematics are described in the total Lagrangian formulation, in a primal-dual approach, including generalized coordinates and generalized velocities. The nonlinear mechanical governing equations are derived based on the variational principle leading to the constrained governing equation for a material body $\setB_0 \subseteq \R^3$,
\begin{equation}
  \label{eq:Hamilton}
  \begin{aligned}
    \delta S
    =
    \int_{\setB_0}
    \bigl( \sprod{\delta\bx}{\dot\bl(\bv,&t)
      + \ff\spint(\bx,t) - \ff\spext(t) + \bH(\bx,t)\tp \blam(t)} \\[-\jot]
    &+ \sprod{\delta\bv}{\bl(\bv,t) - \bl(\dot\bx,t)}
    + \sprod{\delta\blam}{\bh(\bx,t)} \bigr) \, dV
    = 0.
  \end{aligned}
\end{equation}
Here $\sprod{}{}$ denotes a scalar product. $\bx(\theta_i,t)$ and $\bv(\theta_i,t)$ are the positions and velocities depending on the generalized coordinates and generalized velocities of the canonical models, $\blam(t)$ is the vector of Lagrange multipliers of the holonomic constraint $\bh(\bx,t) = 0$ and $\delta(\fcdot)$ denotes their admissible variations. A material point is defined by the coordinates $\theta_i$, with $i = 1,2,3$, in the body's reference coordinate system. The internal and external force densities are denoted by $\ff\spint(\bx,t)$ and $\ff\spext(t)$. $\bH(\bx,t)$ is the Jacobian of kinematical restrictions. Due to the primal-dual formulation, the continuous governing equation enforces the equality of the displacement-based and velocity-based momentum densities by incorporating $\bl(\bv,t) - \bl(\dot\bx,t)$. Inertial forces are considered in the time derivative of the velocity-based momentum density $\dot\bl(\bv,t)$.
\subsection{Spatial and temporal discretization}
The discrete governing equations of the multibody system Finite Element approach are derived through the spatial approximation of the kinematics, specifically the continuous generalized coordinates and velocities of the flexible bodies, using finite elements and solving the integral \eqref{eq:Hamilton}. The temporal discretization is achieved by employing an implicit time integration scheme, which is based on the average vector field method. This approach ensures the preservation of linear and angular momentum as well as the total energy of the system. The methodology is elaborated in Armero et al.\ \cite{armero_energy-dissipative_2003} and has been integrated into the current framework by Gebhardt et al.\ in \cite{gebhardt_new_2020}. Consequently, this leads to the weak form of the discrete governing equations:
\begin{equation}
  \label{eq:Hamilton_discrete}
  \begin{aligned}
    \delta\~S_{n+1}
    &=
    \bigl\langle \delta\~\bq_{n+\frac{1}{2}},
    \dot{\~\bl}(\~\bs_{n+1},\~\bs_n)
    + \~\ff\spint(\~\bq_{n+1},\~\bq_n,\~\bs_{n+1},\~\bs_n) \\
    &- \~\ff\spext(\~\bq_{n+1},\~\bq_n)
    - \~\ff\spext[nc](\~\bq_{n+1},\~\bq_n)
    + \~\bH(\~\bq_{n+1},\~\bq_n)\tp \~\blam_{n+\frac{1}{2}}
    \bigr\rangle \\
    &+ \Sprod{\delta\~\bs_{n+1}}
    {\~\bl(\~\bs_{n+1},\~\bs_n)
      - \~\bl(\~\bq_{n+1},\~\bq_n)} \\
    &+ \Sprod{\delta\~{\blam}_{n+1}}{\~{\bh}(\~\bq_{n+1},\~\bq_n)}
    = 0,
  \end{aligned}
\end{equation}
in which $n$ denotes the discrete time instant at $t_n$. The spatially discrete quantities are indicated by a tilde accent. In \eqref{eq:Hamilton_discrete}, $\~\ff\spext[nc](\~\bq_{t_{n+1}},\~\bq_{t_n})$ specifies the discrete vector of non-conservative forces. In an aeroelastic framework, this vector consists of aerodynamic forces, which are further elucidated in subsequent sections.
\subsection{Computation of aerodynamic loads using the UVLM}
The aerodynamic loads are computed by the unsteady vortex-lattice method (UVLM) based on the work in \cite{Gebhardt2014, Hente_et_al:2023, Katz_Plotkin_2001, roccia2020, Verstraete_2023}. We assume that the influence of viscosity is restricted to a very thin boundary layer. This boundary layer is dominated by vortices which can only be generated there. We further assume the fluid outside of the boundary layer to be incompressible and inviscid, also known as an ideal fluid \cite{Karamcheti:1966}. For flows with a high Reynolds number and a low Mach number, this is a reasonable approximation, which highly simplifies the computation. The fluid behavior can then be described by the well-known Euler equation and the continuity equation for incompressible flows given by
\begin{align}
  \label{eq:Euler}
  \partial_t \bv(\br, t) + (\bv(\br, t) \cdot \nabla) \bv(\br, t)
  &= -\frac1\rho_F \nabla p(\br, t), \\
  \nabla \cdot  \bv(\br, t) &= 0,
  \label{eq:continuity}
\end{align}
where $\bv$ is the fluid particle velocity at its spatial position $\br$ and time instant $t$, $\rho_F$ the fluid density and $p$ the pressure field. Additionally, the velocity field is characterized by a potential flow and can be expressed by Helmholtz's decomposition as the superposition of a scalar potential $\varphi(\br,t)$ and a vector potential $\bm\Psi(\br,t)$ \cite{Bhatia2012} in the following form:
\begin{equation}
  \bv(\br,t)
  =
  \nabla\varphi(\br,t) + \nabla \x \bm\Psi(\br,t)
  =
  \bv_\varphi(\br,t) + \bv_\psi(\br,t),
  \label{eq:helmholtz}
\end{equation}
where the scalar potential component of the velocity is irrotational, and the vector potential component captures any vorticity effect. Furthermore, the kinematical description of the flow can be described by the vorticity field $\bOmega$ defined by
\begin{equation}
  \label{eq:vorticity-definition}
  \bOmega(\br, t) = \nabla \x \bv(\br, t).
\end{equation}
Since flow conditions here are assumed to be dominated by vortices, the fluid velocity is
\begin{equation}
  \label{eq:velocity}
   \bv(\br, t) = \bv_\infty(\br, t) + \bv_B(\br, t) + \bv_W(\br, t),
\end{equation}
with the external flow field velocity $\bv_\infty$, the velocity induced by bounded vortex sheets $\bv_B$ (\eg, boundary of the fluid defined by a structural surface), the velocity induced by free vortex sheets $\bv_W$, \eg, the wake. Without compromising generality, we assume that $\bv_B$ and the free-stream component are absorbed by $\nabla\varphi$ while the field $\bv_W$ is identified with $\nabla\x\bm\Psi$. Incorporating \eqref{eq:helmholtz} into \eqref{eq:continuity} and \eqref{eq:vorticity-definition} facilitates the derivation of a Laplace equation that delineates the fluid's boundary conditions, and a Poisson equation that correlates the potential flow with vorticities \cite{Hente_et_al:2023,Verstraete_2023}. Consequently, the velocities induced by vortices are characterized by the analytical solution provided by the Biot-Savart law. For a vortex segment, this law is defined as follows:
\begin{equation}
  \label{eq:BiotSavarLaw}
  \bv_{\Gamma}(\br, t)
  =
  \frac{\Gamma}{4 \pi} \int_{C(s, t)}
  \frac{\hat{\bm T}(s, t) \x (\br - \br_0(s))}{\norm{\br - \br_0(s)}_2^3}
  \, ds(t),
\end{equation}
where $\hat{\bm T}$ is the unit tangent vector to the curve $C$ and $\Gamma$ is the segment's vorticity.

In the context of the UVLM, the infinitely thin boundary layer around lifting surfaces, where vortices and thus vorticity can be generated or destroyed, is idealized as vortex sheets which in turn are discretized as vortex-lattices with vortex rings consisting of vortex segments having finite circulations. Within each ring, we use vortex segments of equal and constant strength. The computation of the flow field essentially involves the application of the non-penetration boundary condition at the fluid's boundary, specifically at the structural surface. This condition dictates that the normal component $\bm n(\br,t)$ of the relative velocity at the surface should be zero, \ie,
\begin{equation}
  \label{eq:non-Penetration}
  \bv(\br,t) \cdot \bm n(\br,t)
  =
  \bigl[ \bv_\infty + \bv_B(\br, t) + \bv_W(\br, t) - \bv_S(\br, t) \bigr]
  \cdot \bm n(\br,t) = 0,
\end{equation}
where $\bv_S$ is the velocity of the moving surface. Equation \eqref{eq:non-Penetration} reduces the problem to solving for the unknown circulations of the bounded vortex rings. The aerodynamic loads on the surface are ascertained based on \eqref{eq:Euler} solving numerically the unsteady Bernoulli equation, which links the pressure differences across the vortex elements to the prevailing velocity field.
\subsection{Aeroelastic solution}
An important aspect in combining structural and aerodynamic models numerically is the strategy adopted for the transfer of information between their respective meshes.  We transfer the aerodynamic loads coming from the UVLM into the structural model, stating that for any time the virtual work done by the aerodynamic loads on the aerodynamic mesh is equal to the virtual work done of the generalized aerodynamic forces $\~\ff\spext[ae]$ on the nodes of the structural mesh \cite{Gebhardt2014, Hente_et_al:2023}. The spatial coordinates of any point on the fluid domain are mapped into the configuration space of the structural model using weighted linear surjective vector-valued mapping functions comprising radial-based $C^\infty$ bump functions with compact support \cite{Hente_et_al:2023}.

Following \eqref{eq:Hamilton_discrete}, the discrete nonlinear governing equation in a strong  fluid-structure interaction, wherein the structural behavior and aerodynamics are concurrently solved at an identical time instant in terms of the structural variables, can now be stated as follows:
\begin{equation}
  \ff(\~\bq, \~\bs, \~\blam)_{n+1}
  =
  \col[c]{
    \~\ff\spint(\~\bq,\~\bs) - \~\ff\spext[ae](\~\bq,\~\bs) +
    \dot{\~\bl}(\~\bs) + \~\bH(\~\bq)\tp \~\blam \\
    \~\bl(\~\bs) - \~\bl(\~\bq) \\
    \~\bh(\~\bq)
  }_{n+1}
  =\bm0,
  \label{eq:system-str-aero}
\end{equation}
where the subscript $n+1$ indicates that the unknowns are solved at time step $t_{n+1}$. We do not explicitly express \eqref{eq:system-str-aero} with the specific time dependencies arising from the employed time integration method, denoted as $(\~\bq_{n},\~\bq_{n+1}, \~\bs_{n},\~\bs_{n+1},\~\blam_{n+1/2})$ and outlined in \eqref{eq:Hamilton_discrete}. However, it is important for the reader to bear in mind these temporal dependencies when considering the system's dynamics. In \eqref{eq:system-str-aero}, the first equation represents the discrete dynamic equilibrium, capturing the balance of forces acting on the system. The second equation is the momentum equivalency, and the third one represents the discrete constraints equations. We solve the nonlinear governing equation iteratively applying Newton's method, deriving its linearized form with the Taylor expansion by neglecting the higher order terms to:
\begin{equation}
  \ff(\~\bq, \~\bs, \~\blam)_{n+1}^{k+1}
  =
  \ff(\~\bq, \~\bs, \~\blam)_{n+1}^{k} +
  \Delta \ff(\~\bq, \~\bs, \~\blam)_{n+1}^{k} = \bm{0},
  \label{eq_ae_newton_to_governing}
\end{equation}
where $k$ stands for the iteration step within the Newton iteration process and $\Delta \ff$ is the discrete increment of $\ff$ obtained by calculating the partial derivatives with respect to the discrete generalized variables and the Lagrange multiplier, \ie,
\begin{align}
  \Delta \ff(\~\bq, \~\bs, \~\blam) _{n+1}^{k}
  &= \pfrac{\ff}{\~\bq} \Delta \~\bq +
    \pfrac{\ff}{\~\bs} \Delta \~\bs +
    \pfrac{\ff}{\~\blam} \Delta \~\blam
    =
    \bK(\~\bq, \~\bs. \~\blam)_{n+1}^{k} \Delta\~\by ,
    \label{eq_ae_newton_incremental}
\end{align}
In \eqref{eq_ae_newton_incremental},
$\~\bK$ is the derivative matrix (or Jacobian) of \eqref{eq:system-str-aero}
and $\Delta \~\by=(\Delta\~\bq, \Delta\~\bs, \Delta\~\blam)$.
We evaluate all partial derivatives in the linear system
\eqref{eq_ae_newton_incremental} analytically as described in
\cite{hente_modal_2019} and \cite{Hente_et_al:2023}.
We solve the system for all variables, including the Lagrange multiplier.
Within the aeroelastic approach, the matrix consists of two parts:
firstly, the Jacobian of the structural model,
derived from the partial derivatives of the discrete structural model,
the discrete equivalency of linear momentum, and the constraint equations;
and secondly, the Jacobian of the discrete generalized aerodynamic forces,
which involves the partial derivatives of the aerodynamic forces.
For an exhaustive explanation and derivation of both matrices,
readers are referred to
\cite{Gebhardt_et_al:2020b,hente_modal_2019,hente_static_2021,Hente_et_al:2023}.


\section{Solution Algorithms}
\label{sec:solution}
In this section we focus on solving the nonlinear equation system
for an integration step of the strongly coupled aeroelastic system
within an implicit integration scheme.
We start with a brief recapitulation of the relevant Newton techniques as such
before discussing their application in our aeroelastic simulation.

Given a smooth function $\bF\: \R^n \to \R^n$,
we seek a solution $\by^*$ of the nonlinear equation system
\begin{equation}
  \label{eq:nonlinear-system}
  \bF(\by) = 0.
\end{equation}
The standard Newton iteration starts with an initial guess $\by_0$
and updates each iterate $\by_k$ with a step $\incN$
that solves the linearization of \eqref{eq:nonlinear-system} at $\by_k$,
called the Newton system,
\begin{equation}
  \label{eq:newton-system}
  \by_{k + 1} = \by_k + \incN
  \qtextq{where}
  \bF'(\by_k) \incN = -\bF(\by_k).
\end{equation}
The key property of Newton's method is quadratic local convergence,
which is guaranteed if the derivative matrix $\bF'(\by) \in \R^{n \x n}$ is
invertible and $\bF'$ is Lipschitz continuous in a convex domain around $\by^*$.
Modern affine-invariant convergence theory even provides sharp estimates
of the convergence speed and convergence radius in terms of the
Lip\-schitz constant (which measures the maximal curvature of $\bF$)
and of the initial step $\incN[0]$ \cite{Deuflhard:2004};
see also \cite{Bock:1987} for the seamless extension
to constrained nonlinear least-squares problems.
Moreover, global convergence can be enforced
with a steplength control,
\begin{equation}
  \label{eq:newton-damped}
  \by_{k + 1} = \by_k + \alpha_k \incN,
  \quad
  \alpha_k \in (0, 1],
\end{equation}
where the steplength $\alpha_k$ is determined by a standard
backtracking line search based on the Armijo condition,
\begin{equation}
  \label{eq:newton-armijo}
  \norm{\bF(\by_k + \alpha_k \incN)}_2^2
  \le
  (1 - \tfrac{c}{2}) \norm{\bF(\by_k)}_2^2
  \qtextq{with} c \in (0, 1).
\end{equation}
This \emph{damped} Newton iteration converges to a
stationary point $\by^*$ of $\norm{\bF(\by)}_2$ under mild assumptions.
In the following we refer to the standard Newton method as
\emph{exact Newton method} and reserve the notation $\incN$
for the exact Newton step.

Another essential property of Newton's method is that it suffices
to solve the linear system \eqref{eq:newton-system} approximately:
small errors in $\incN$ do not destroy the local convergence,
although the convergence rate will usually be
superlinear or merely linear instead of quadratic.
An easy way to exploit this property consists in replacing
each matrix $\bF'(\by_k)$ with an invertible approximation~$\bB_k$,
either an ad hoc choice or a systematic choice,
when the exact derivative is expensive to evaluate or to factorize.
Quasi-Newton methods represent the most prominent
systematic approach in this spirit.
Originally proposed by Broyden \cite{Broyden:1965}
based on prior work of Davidon \cite{Davidon:1959},
the key idea is to construct $\bB_{k+1}$ from $\bB_k$
by some cheap update formula such that the \emph{secant condition} holds,
\begin{equation}
  \label{eq:secant}
  \bB_{k+1} \incr = \bF(\by_{k+1}) - \bF(\by_k),
\end{equation}
where the quasi-Newton step $\incr \approx \incN$ solves
\begin{equation}
  \label{eq:qn-system}
  \bB_k \incr = -\bF(\by_k).
\end{equation}
The matrix $\bB_{k+1}$ is not uniquely determined by \eqref{eq:secant},
and different additional requirements lead to different quasi-Newton methods.
Moreover, it is possible to work directly with the inverse $\bH_k = \Inv \bB_k$
to obtain the step $\incr = -\bH_k \bF(\by_k)$ without solving a linear system.
These quasi-Newton methods are primarily employed
for minimizing (or maximizing) a function $\phi\: \R^n \to \R$
by solving the stationarity condition $\phi'(\by) = 0$.
Here we have $\bF(\by) = \phi'(\by)$, the Jacobian of $\bF$ is the
\emph{symmetric} Hessian of $\phi$, $\bF'(\by) = \phi''(\by)$,
and the convergence theory is even richer than for general equations
\cite{
Dennis_Schnabel:1983,Dennis_More:1977,Fletcher:1987}.
In the sequel we will use the term \emph{quasi-Newton method}
for any Newton type iteration
that replaces $\bF'(\by_k)$ in \eqref{eq:newton-system}
with some approximation $\bB_k$
to compute the step $\incr$ from \eqref{eq:qn-system}.

Of course, approximate solutions of \eqref{eq:newton-system}
can also be obtained by iterative methods,
which leads to the class of \emph{inexact Newton methods}.
Originally proposed by Dembo, Eisenstat and Steihaug
\cite{Dembo_et_al:1982,Dembo_Steihaug:1983},
these methods have been further developed and analyzed by Deuflhard
\cite{Deuflhard:1991,Deuflhard:2004} and others.
They produce a finite sequence of increasingly accurate approximations
$\incr^1, \dots, \incr^i =:\incr$
to the current exact step $\incN$, thus offering
the advantage of controlling the residual error $\br_k$
of the inexact Newton step $\incr \approx \incN$,
\begin{equation}
  \label{eq:inexact}
  \br_k = \bF'(\by_k) \incr + \bF(\by_k),
\end{equation}
such that it matches the progress of the Newton iteration,
\begin{equation}
  \label{eq:forcing}
  \norm{\br_k} \le \eta_k \norm{\bF(\by_k)},
  \quad
  \eta_k \in (0, 1).
\end{equation}
A proper adaptive choice of the \emph{forcing sequence} $(\eta_k)_{k \in \N}$
then guarantees superlinear or even quadratic local convergence.
Thus inexact Newton methods can reduce the effort
for solving the linear system \eqref{eq:newton-system}.
However, unlike quasi-Newton methods,
they require the evaluation of the derivative $\bF'(\by_k)$
or at least of the product $\bF'(\by_k) \incr$ at every iterate.

In our setting, the quasi-Newton and inexact Newton algorithms
use the same approximate Jacobian, $\bB_k \approx \bF'(\by_k)$,
and the trial steps $\incr^j$ are generated by iterative refinement,
\begin{equation*}
  \incr^j = \incr^{j-1} - \Inv\bB_k \br_k^{j-1},
  \quad
  \br_k^j = \bF'(\by_k) \incr^j + \bF(\by_k),
\end{equation*}
starting with $\incr^0 = 0$ and hence $\br_k^0 = \bF(\by_k)$.
Local convergence of the inexact Newton approach is then guaranteed
if the approximation quality of each $\bB_k$ is sufficiently good
in the rigorous sense that $\bB_k$ gives a spectral radius
\begin{equation*}
  \rho_k := \rho(\Inv\bB_k \bF'(y_k) - I) \stackrel!< 1,
\end{equation*}
because this implies convergence of $\incr^j$ to $\incN$,
\begin{equation*}
  \norm{\br_k^j} \le
  \rho_k \norm{\br_k^{j-1}} \le
  \rho_k^j \norm{\bF(\by_k)} \to 0.
\end{equation*}
The quasi-Newton approach, where $\incr = \incr^1$,
can still converge slowly in this situation
because it does not control the reduction of the residual error,
and it can even diverge if
$\norm{\bF(\by_k + \incr^1)} \ge \norm{\bF(\by_k)}$.
The inexact Newton approach is therefore preferable
from a theoretical perspective.

We are now ready to proceed with the nonlinear
Newton system \eqref{eq:system-str-aero}
of the implicit time step for the aeroelastic simulation.
Omitting the tilde for simplicity,
the (discrete) variables are $\by = (\bq, \bs, \blam)$.
The derivative $\bF'(\by)$ consists of two parts:
one from the mechanical structure and
one from the aerodynamics,
\begin{align}
  \bF'(\by)
  &= \Mstr(\by) + \Maer(\by)
  \notag \\
  \label{eq:str+aer}
  &=
  \mat[3]{
    \bK_{\bq\bq}(\by) & \bK_{\bq\bs}(\by) & \bH_d(q)\tp \\
    \bK_{\bs\bq}(\by) & \bK_{\bs\bs}(\by) & 0 \\
    \bH_d(q) & 0 & 0
  }
  +
  \mat[3]{
    \bK\tsp{aero}_{\bq\bq}(\by) & \bK\tsp{aero}_{\bq\bs}(\by) & 0 \\
    0 & 0 & 0 \\
    0 & 0 & 0
  }
  .
\end{align}
The structural derivative $\Mstr(\by)$ possesses a near-symmetric
saddle-point structure with sparse blocks.
It is comparatively cheap to evaluate and to factorize.
The aerodynamic derivative $\Maer(\by)$ is non-symmetric
with two nearly dense blocks,
which are very expensive to evaluate and
which have small entries as compared to $\Mstr(\by)$.
Moreover, the full derivative $\bF'(\by)$
is considerably more expensive to factorize
than just $\Mstr(\by)$ because of the two dense blocks.
The individual matrix blocks are defined as follows.
The first row contains partial derivatives of the dynamic equilibrium,
giving structural blocks
\begin{align*}
  \bK_{\bq\bq}(\by) &= \pfrac{\ff\spint(\bq,\bs)}{\bq}, &
  \bK_{\bq\bs}(\by) &= \pfrac{[\ff\spint(\bq,\bs) + \dot\bl(\bs)]}{\bs},
\end{align*}
and aerodynamic blocks
\begin{align*}
  \bK\tsp{aero}_{\bq\bq}(\by) &= -\pfrac{\ff\spext[ae](\bq,\bs)}{\bq}, &
  \bK\tsp{aero}_{\bq\bs}(\by) &= -\pfrac{\ff\spext[ae](\bq,\bs)}{\bs}.
\end{align*}
The second row contains partial derivatives of the momentum equivalence, giving
\begin{align*}
  \bK_{\bs\bq}(\by) &= -\bl'(\bq), & \bK_{\bs\bs}(\by) &= \bl'(\bs).
\end{align*}
Finally the Jacobian of the holonomic constraint,
$\bH_d(\bq) = \bh'(\bq)$, appears in the lower left
and (transposed) in the upper right of $\Mstr(y)$.

We consider three solution approaches for the nonlinear equation system:
\begin{enumerate}
\item the exact Newton method wherein \eqref{eq:newton-system}
  is solved directly;
\item a quasi-Newton method wherein \eqref{eq:qn-system}
  is solved directly with $\bB_k = \Mstr(\by_k)$ \cite{Hente_et_al:2023};
\item an inexact Newton method wherein \eqref{eq:newton-system}
  is solved by iterative refinement
  using the structural inverse $\Inv{\Mstr(\by_k)}$
  as preconditioner for the full inverse $\Inv{\bF'(\by_k)}$.
\end{enumerate}
At each iteration,
approach (1) requires one factorization of $\bF'(\by_k)$ and one backsolve,
approach (2) requires one factorization of $\Mstr(\by_k)$ and one backsolve, and
approach (3) requires one factorization of $\Mstr(\by_k)$ and several backsolves.
In all cases we employ Pardiso \cite{Schenk_Gaertner:2011}
from the Intel Math Kernel Library (MKL)
as direct sparse solver for the factorization and the backsolves.

The cost per iteration is very high for approach (1)
since the entire partially dense matrix $\bF'(\by_k)$
needs to be evaluated and factorized.
For approach (2) the cost per iteration is very low
since the expensive dense blocks of $\Maer(\by)$ are never evaluated
and only the sparse matrix $\Mstr(\by_k)$ needs to be factorized.
For approach (3) the cost per iteration is in between
since again only $\Mstr(\by_k)$ needs to be factorized
but $\Maer(\by)$ needs to be evaluated for the refinement steps.


\section{Computational Results}
\label{sec:results}
We consider two computational examples:
a thin, flexible plate in an external airflow
and a wind turbine in operating conditions.
All computations are performed on a single core of a desktop PC
with eight AMD Ryzen 7 7700 cores and \SI{64}{GiB} of RAM.
The simulation code \DeSiO is implemented in Fortran 90,
the inexact Newton algorithm is implemented in \CC,
and the time measurements are performed using the
\CC library \sfname{boost::chrono} \cite{Boost}.

The time steps in both examples are so small that each Newton iteration
starts in the local convergence domain;
therefore we perform full-step iterations without any steplength control.
Moreover, the employed approximate inverse $\Inv{\Mstr(\by)}$
is quite close to the full inverse $\Inv{\bF'(\by)}$ throughout all simulations.
This means that all three Newton algorithms converge safely and rapidly,
which allows for a direct comparison of computation times.


\subsection{Plate}
\label{sec:plate}
\begin{figure}
  \centering
  \ximg[width=0.9\linewidth,trim=32 10 0 15,clip]{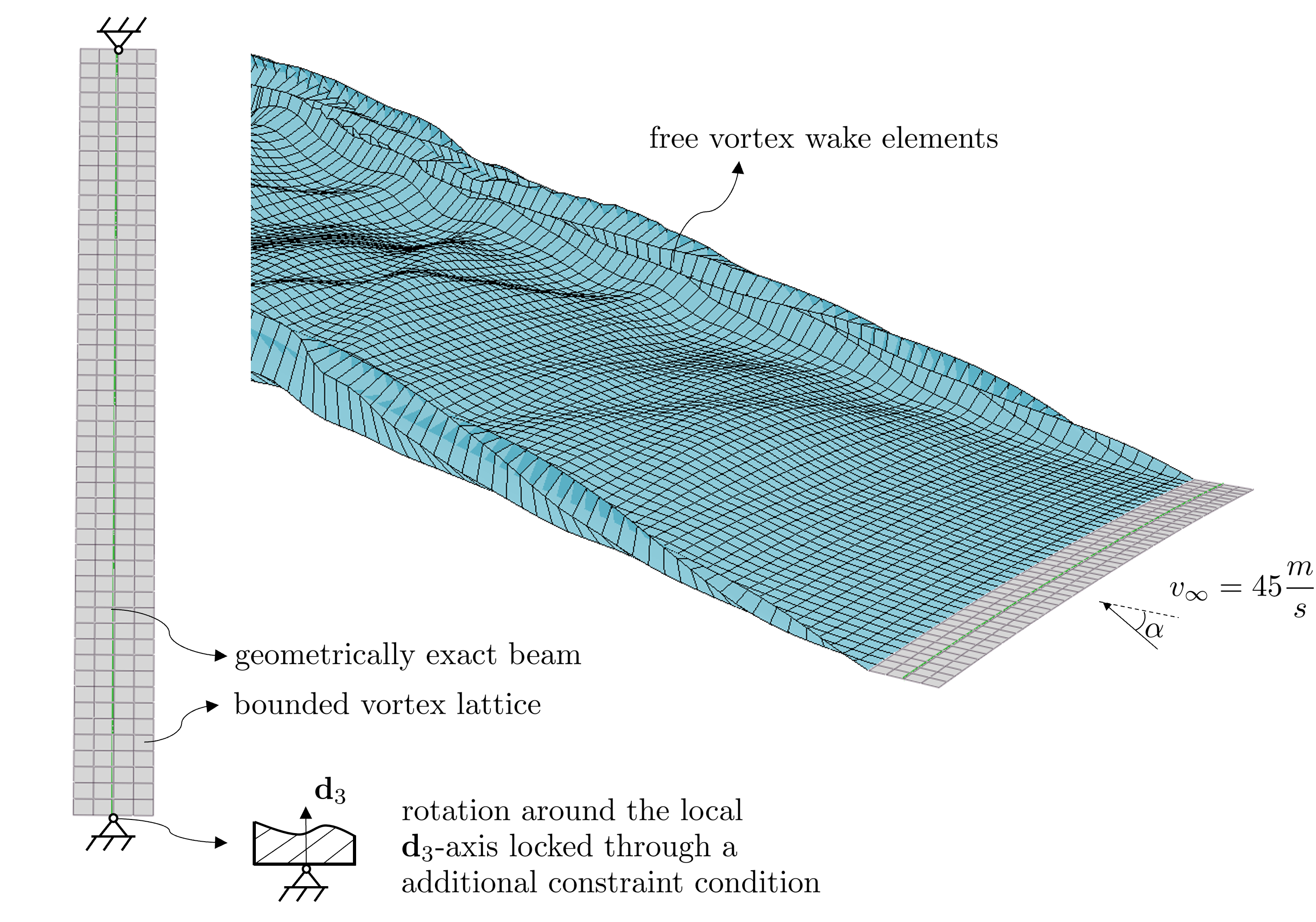}
  \caption{Illustration of the rigidly supported flexible plate.}
  \label{fig:plate}
\end{figure}
The following analysis explores the dynamic aeroelastic response of a
thin, flexible aluminum plate when subjected to an external airflow.
The plate, measuring $\SI{10.0}{m}$ in length,
$\SI{1.0}{m}$ in width,  and $\SI{8.0}{mm}$ in thickness,
adheres to standard isotropic linear elastic properties,
featuring a Young's modulus of $\SI{7.0e10}{N/m^2}$,
a shear modulus of $\SI{2.63e10}{N/m^2}$,
and a material density of $\SI{2.7e3}{kg/m^3}$.
The plate's free ends are hinged supported,
and the torsional degree-of-freedom at the left end is constrained as well.
The air inflow is characterized by a free-stream velocity,
denoted as $V_\oo = \SI{45}{m/s}$,
and an angle of attack $\alpha = \SI{15}{\degree}$.
The flow is applied impulsively
and the air density is specified as $\rho_F = \SI{1.225}{kg/m^3}$.
A summary of the geometric and material input data
for both the structural and aerodynamic models
is given in \cref{tab:ex01_main_geometry}.
\begin{table}
  \centering
  \caption{Geometry and material parameters for beam model of flexible plate.}
  \label{tab:ex01_main_geometry}
  \begin{tabular}{lr}
    \toprule
    Parameter & Value \\
    \midrule
    Length $L$    &  \SI{10.0}{m} \\
    Chord $c$     &   \SI{1.0}{m} \\
    Thickness $t$ & \SI{0.008}{m} \\[\jot]
    Young's modulus $E$     &  \SI{7.0e10}{N/m^2} \\
    Shear modulus $G$	    & \SI{2.63e10}{N/m^2} \\
    Material density $\rho$ &  \SI{2.70e3}{kg/m^3} \\
    Fluid density $\rho_F$  &   \SI{1.225}{kg/m^3} \\
    \bottomrule
  \end{tabular}
\end{table}
In our aeroelastic model, we discretized the infinitely thin boundary layer
using a mesh of $m_A \x n_A = 50 \x 4$ vortex-rings,
where $m_A$ denotes the span-wise discretization
and $n_A$ represents the chord-wise discretization.
The structural model is characterized by $m_S = 50$
geometrically exact beam elements.
As mentioned above, the information transfer between both models
is ensured by a radial-based function.
For this, it is necessary to define a fixed transferring radius
which serves as a distance threshold,
limiting the range of interest in transferring information among the models.
For this example, we use a search radius of $\gamma\tsb{ref} = \SI{0.501}{m}$
at each node.

The simulation is conducted until a steady-state solution is attained,
necessitating a total simulation time of $T_n = \SI{3}{s}$.
The convection of the free wake is achieved
considering a characteristic length $\Delta L = \SI{0.25}{m}$
resulting in a time increment $\Delta t = \SI{5.56e-3}{s}$.
To mitigate the effects of singularities
resulting from vortex-induced velocities,
we introduce a cut-off value, denoted as $\delta = 0.01$.
The main simulation parameters are summarized in
\cref{tab:ex01_model_parameter}.
\begin{table}
  \centering
  \caption{Simulation parameters for aeroelastic analysis of flexible plate.}
  \label{tab:ex01_model_parameter}
  \begin{tabular}{lr}
    \toprule
    Parameter & Value \\
    \midrule
    Total simulation time $T_n$             & \SI{3}{s} \\
    Simulation time increment $\Delta t$    & \SI{5.56e-3}{s} \\
    Cut-off radius $\delta$                 & \num{0.01} \\
    Transfer-radius $\gamma\tsb{ref}$       & \SI{0.501}{m} \\
    Angle of attack $\alpha$                & \SI{15.0}{\degree} \\
    Intensity of free-field velocity $V_\oo$ & \SI{45}{m/s} \\
    Newton convergence tolerance tol$_N$    & \num{1e-8} \\
    \bottomrule
  \end{tabular}
\end{table}

\begin{figure*}[tp]
  \centering
  \fboxsep=0pt
  \fboxrule=0pt
  \begin{subfigure}[b]{0.6\textwidth}
    \centering
    \fbox{\includegraphics[width=\linewidth,trim=60 4 120 36,clip]
      {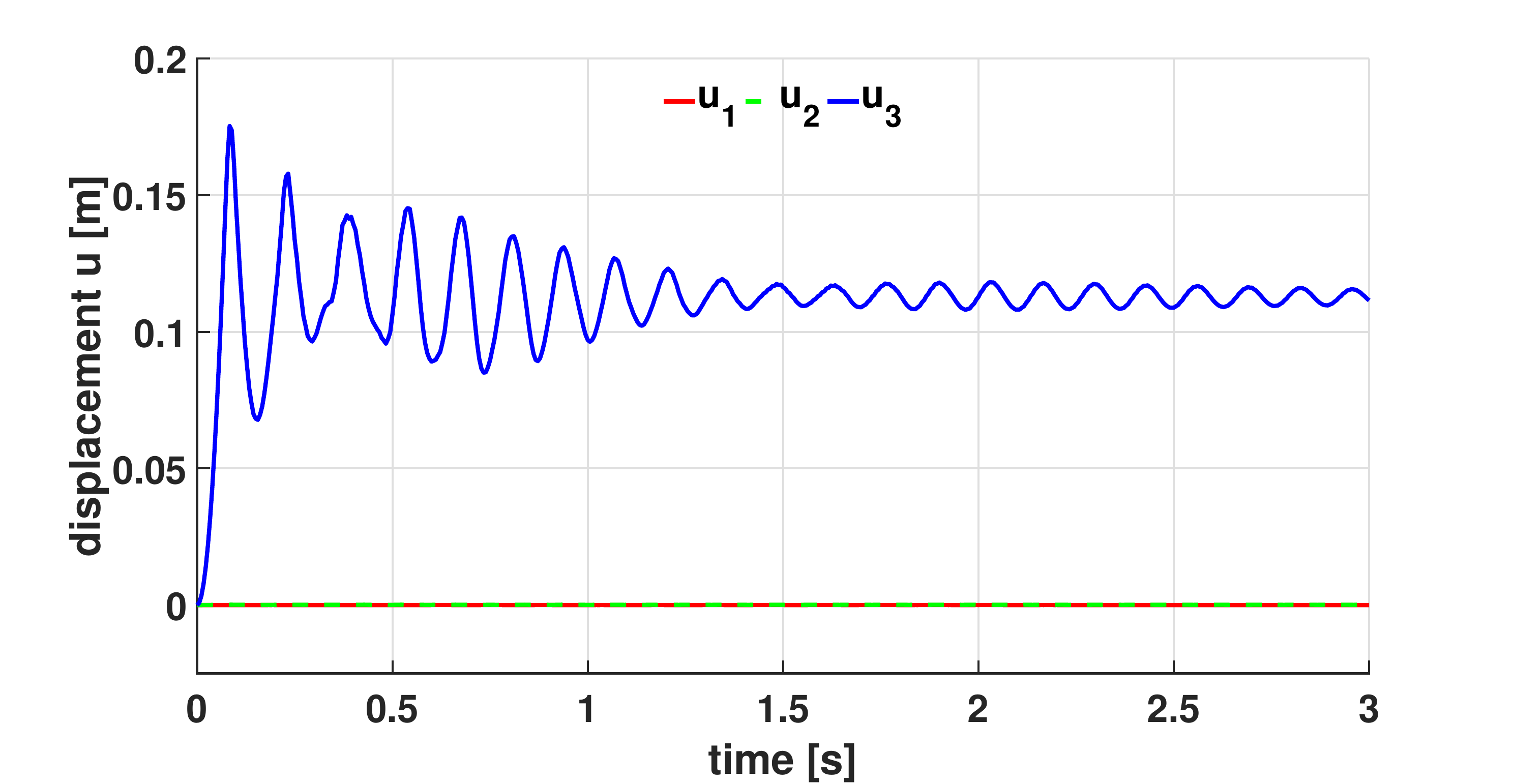}}
    \caption{Displacement of the beam center over time.}
    \label{fig:plate_results-displacement}
  \end{subfigure}
  \hfill
  \begin{subfigure}[b]{0.6\textwidth}
    \centering
    \fbox{\includegraphics[width=\textwidth,trim=60 4 120 36,clip]
      {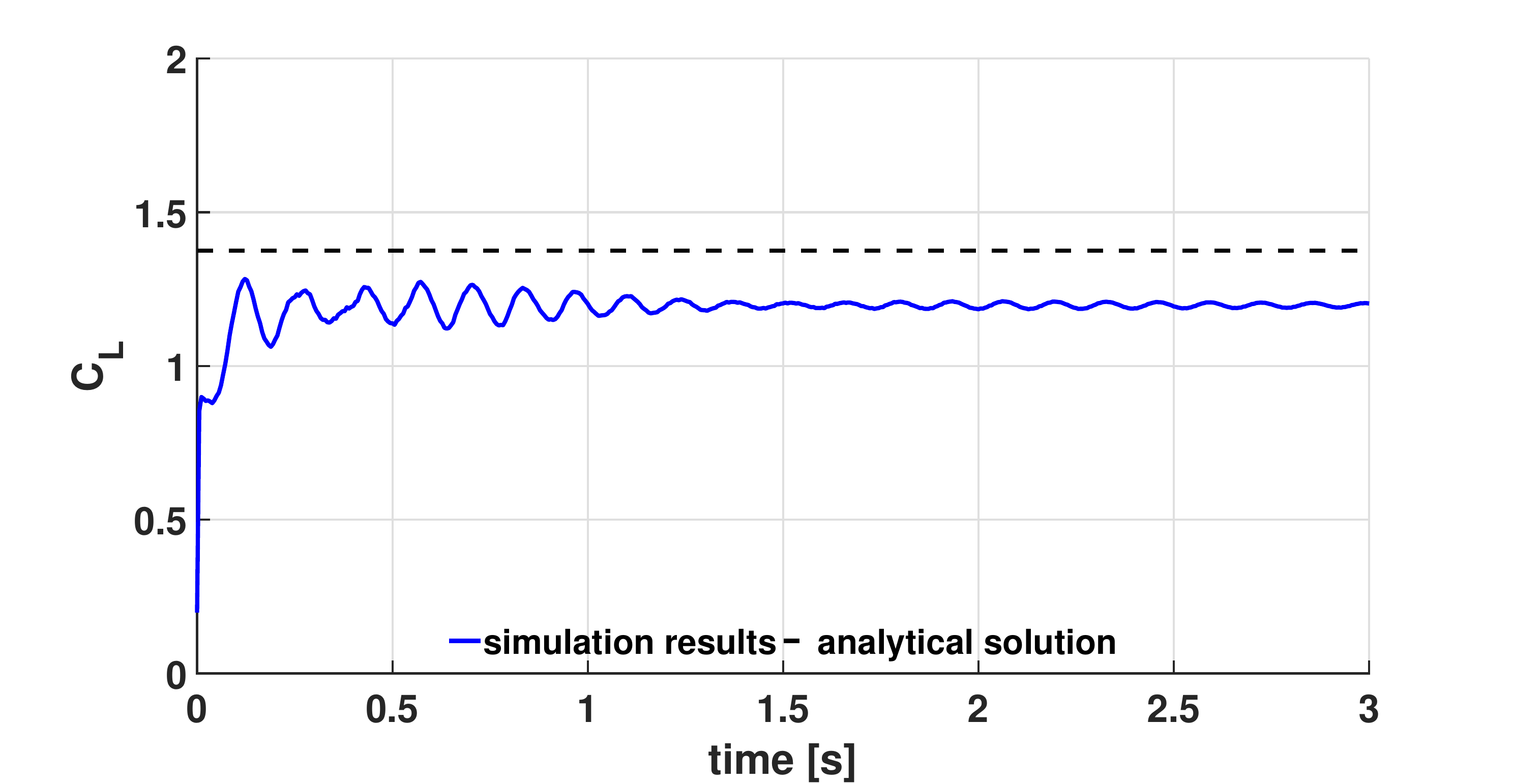}}
    \caption{Lift coefficient over time.}
    \label{fig:plate_results-CL}
  \end{subfigure}
  \caption{Results of the aeroelastic simulation of the flexible plate.}
  \label{fig:plate_results}
\end{figure*}

In \cref{fig:plate_results}, we show some physical results
for validation purposes.
On top, in \cref{fig:plate_results-displacement},
the displacement of the node at the center of the beam over time is given.
It can be seen that the aerodynamic forces first excite the structure
before aerodynamic damping ensures that a steady state is approached.
The displacement occurs mainly in the vertical direction,
i.e., almost no rotation takes place.
A similar behavior is observed with regard to the lift coefficient $C_L$,
shown at the bottom in \cref{fig:plate_results-CL}.
After some initial vibrations, the system approaches a steady state.
The resulting angle of attack is shown together with the analytical solution
for a finite plate \cite{Katz_Plotkin_2001}
for the resulting steady-state angle of attack, $\alpha = 15.05^\circ$.
It should be noted that this analytical solution is derived
for elliptic airfoils and small angles of attack.
Therefore, the small deviation of ca.\ 12.5\%
in the resulting lift coefficient is acceptable.
It should also be mentioned that the settings for our simulation
are chosen not to ensure a description as physically accurate as possible
but to test the different solution algorithms.
As expected, the results for all chosen solution algorithms are identical
since the same system of equation is solved.

The absolute tolerance for all three Newton algorithms is $10^{-8}$,
\ie, each iteration terminates when $\norm{\bF(\by_k)} \le 10^{-8}$.
For the inexact Newton algorithm we consider two different forcing sequences,
\begin{align*}
  \eta_k^1 &= \min\Bigl(\frac12, \norm{\bF(\by_k)}\Bigr),
  &\eta_k^2 &= 10^{-5} \eta_k^1.
\end{align*}
Both choices guarantee quadratic convergence.
The computational results are presented in \cref{tab:plate:solvers},
where ``Inexact Newton $i$'' refers to forcing sequence $(\eta_k^i)$.
\begin{table}
  \let\mc\multicolumn
  \def\N{\rlap.\enspace\ignorespaces}
  \def\M{\N\N}
  \centering
  \caption{Flexible plate: runtimes in seconds and iteration counts.
    Inexact Newton $i$ refers to forcing sequence $(\eta_k^i)$.
    Integration = complete simulation (540 time steps),
    Newton = solving nonlinear system \eqref{eq:system-str-aero},
    Eval UVLM = evaluating aerodynamic data ($\Maer$ and $\faer$ in Newton),
    Eval structure = evaluating mechanical data
    ($\Mstr$ and $\ff\tsp{str}$ in Newton).}
  \label{tab:plate:solvers}
  \begin{tabular}{l*4r}
    \toprule
    & Exact  & Quasi- & Inexact  & Inexact  \\
    & Newton & Newton & Newton 1 & Newton 2 \\
    \midrule
    \mc5c{\it Total} \\[1pt]
    Integration       & 4135.22 & 2362.30 & 4314.52 & 3997.82 \\
    \N Eval UVLM      & 2094.09 & 2099.47 & 2113.71 & 2103.91 \\
    \N Eval structure &    0.23 &    0.22 &    0.23 &    0.23 \\
    \N Newton         & 2040.89 &  262.60 & 2200.57 & 1893.68 \\
    \M Eval UVLM      & 1872.54 &  246.61 & 2181.24 & 1874.14 \\
    \M Eval structure &    2.84 &    3.09 &    3.29 &    2.83 \\
    \M Linear solver  &  163.53 &   12.75 &   15.98 &   16.66 \\[2pt]
    Newton steps      &    2625 &    3041 &    3038 &    2625 \\
    Refinement steps  &    ---  &     --- &    4375 &    8594 \\
    \midrule
    \mc5c{\it Average per time step} \\[2pt]
    Integration       &  7.657 &  4.374 &  7.989 &  7.403 \\
    \N Eval UVLM      &  3.878 &  3.888 &  3.914 &  3.896 \\
    \N Newton         &  3.779 &  0.486 &  4.075 &  3.506 \\
    \M Eval UVLM      &  3.467 &  0.456 &  4.039 &  3.470 \\
    \M Eval structure & .00524 & .00572 & .00607 & .00522 \\
    \M Linear solver  & .30281 & .02361 & .02959 & .03085 \\[2pt]
    Newton steps      &  4.861 &  5.631 &  5.625 &  4.861 \\
    Refinement steps  &   ---  &    --- &  8.102 & 15.915 \\
    \bottomrule
  \end{tabular}
\end{table}

We observe that in all cases at most half the total simulation time
is spent on solving the nonlinear equation system \eqref{eq:system-str-aero}
by the respective Newton algorithms. The remaining computing time
(roughly \SI{2100}{s} total or \SI{3.9}{s} per time step)
is almost entirely spent on UVLM evaluations during the time integration.
Within the Newton algorithms the time is also predominantly spent
on UVLM evaluations, \ie, on arodynamic forces and,
except for quasi-Newton, on the associated derivative matrices.
We focus on the computation time spent in the Newton algorithms.
Herein, the actual solving (of linear systems)
amounts to less than $10$\% (exact), $5$\% (quasi-Newton),
and $1$\% (inexact), respectively.

The quasi-Newton algorithm is clearly
the fastest nonlinear solver in this simulation.
On the one hand, it generates only approximate Newton increments,
which increases the number of required Newton steps by roughly $20$\%
as compared to the exact Newton algorithm.
On the other hand, it avoids the evaluation
of aerodynamic derivatives $\Maer(\by_k)$ entirely
and in addition it solves the modified linear system much faster,
which ultimately reduces the total solution time to just $13$\%.
There is only one drawback:
the performance of the quasi-Newton method depends entirely
on the quality of all encountered approximate inverses,
which is not controlled or even measured;
therefore the algorithm can become very slow or even diverge
in unfortunate cases.

The inexact Newton algorithm avoids this drawback,
and with the first forcing sequence
it often produces the same Newton increments
as the quasi-Newton algorithm,
taking only one or two linear system solutions per Newton step.
However, the additional control of the linear system residual
requires the evaluation of aerodynamic derivatives $\Maer(\by_k)$.
Thus, although the linear system solving
is just moderately slower than with the quasi-Newton algorithm
and still ten times as fast as with the exact Newton algorithm,
the total solution time increases
beyond that of the exact Newton algorithm
because of the additional $20$\% Newton steps.
This is a very peculiar situation:
the time for linear system solving is virtually irrelevant
since the time for evaluationg UVLM data dominates everything.

The previous observations motivate
the choice of the second forcing sequence.
Here we increase the effort spent in the (cheap) linear solver
by asking for five additional decimal digits of accuracy
for every Newton increment.
In effect, this yields essentially the same accuracy
as in the exact Newton algorithm
and precisely the same number of Newton steps: 2625 instead of 3038.
Now the number of iterative refinement steps almost doubles
as compared to the first forcing sequence,
but the time for linear system solution increases by less than $5$\%
while the total solution time decreases by $14$\%,
or by $7$\% as compared to the exact Newton algorithm.

In summary, with a succinct choice of the forcing sequence,
the inexact Newton algorithm outperforms the exact Newton algorithm
while matching its accuracy and guaranteeing quadratic convergence
(provided that the structural inverse $\Inv{\Mstr(\by_k)}$
is sufficiently accurate for convergence of the iterative refinement).
The quasi-Newton algorithm outperforms both of them by far,
although with no convergence guarantee at all.


\subsection{Wind Turbine}
\label{sec:turbine}
As a second example we consider the transient simulation
of a wind turbine in operating conditions.
We use the IEA \SI{15}{MW} reference wind turbine \cite{IEA_2020_15MW}
with a rotor diameter of \SI{240}{m} in this example.
Since we are only interested in the aeroelastic simulation,
we do not consider any hydrodynamic forces or control input,
instead considering the turbine fixed at the bottom of the tower.
To avoid an initial impulse on the system,
we start the aeroelastic simulation in a state
based on the results of two structural pre-simulations.
We first consider the deformation due to gravitity and, based on this,
we then perform the structural computations for the rotor
to ensure that the velocities transferred from the structural model
to the aerodynamic model are correct in the first time step.

The wind speed is constant at $V_\infty = \SI{8}{m/s}$,
and starting from its initial velocity the rotor can rotate freely,
\ie, we do not consider any generator torque.
According to the report \cite{IEA_2020_15MW},
we set the pitch angle to a fixed value of $\varphi = \SI{0}{\degree}$
and the initial rotation speed to \SI{5.6}{rpm}.
These conditions are between the cut-in wind speed
and the rated wind speed of the turbine.

\begin{figure}
  \centering
  \ximg[width=\linewidth,trim=8 32 10 15,clip]{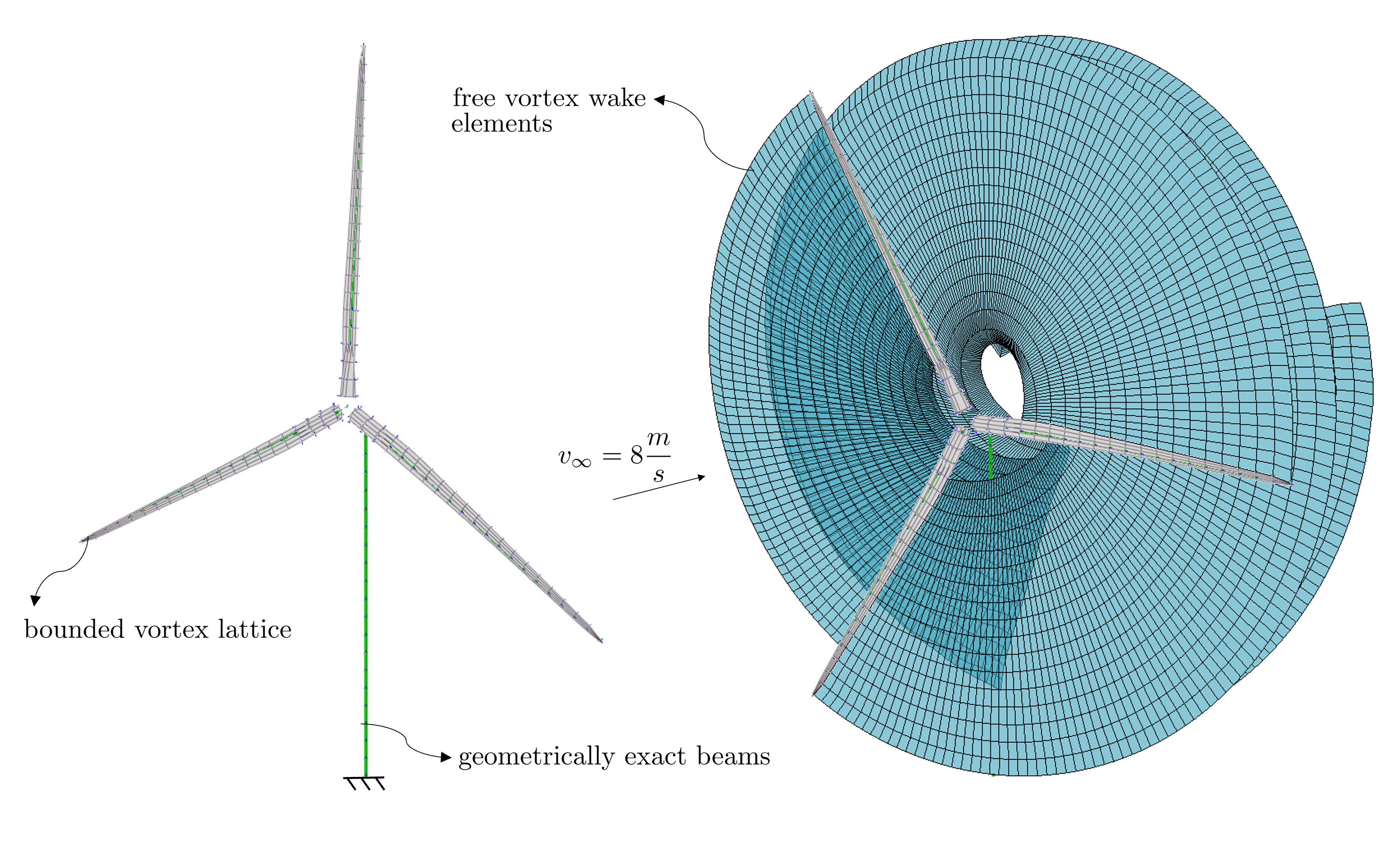}
  \caption{Structural (left) and aeroelastic (right)
    model of the IEA \SI{15}{MW} reference wind turbine
    including wake in DeSiO (blade surfaces shown for visualisation only).}
    \label{fig:desio_modeshape}
\end{figure}

The wind turbine model consists of four geometrically exact beams
for the tower and the three blades
while the hub and nacelle are modelled as rigid bodies.
The material and geometrical properties of all components
are again based on the report \cite{IEA_2020_15MW}.
Only the blades are considered in the computation of the aerodynamic forces;
the influence of tower, hub, and nacelle is neglected.
Each beam consists of 20 elements,
the vortex sheet is discretised by 20 elements in span-wise
and 4 elements in chord-wise direction.
A simulation time of 10 seconds is considered
with a time step $\Delta t = \SI{0.05}{s}$.
The resulting wake is illustrated in \cref{fig:desio_modeshape}.

The absolute tolerance for all three Newton algorithms is now $10^{-6}$,
and for the inexact Newton algorithm we now consider just one forcing sequence
(which again guarantees quadratic convergence):
$\eta_k = \min(\frac12, \norm{\bF(\by_k)})$.
The computational results are presented in \cref{tab:wea:solvers}.
\begin{table}
  \let\mc\multicolumn
  \def\full{\llap{$\to$\kern.06em}100}
  \def\time#1{Time$\hphantom{#1}$}
  \def\N{\rlap.\enspace\ignorespaces}
  \def\M{\N\N}
  \centering
  \caption{Wind turbine: runtimes in seconds and in percent of
    exact Newton (\%E), quasi-Newton (\%Q), inexact Newton (\%I),
    and iteration counts.
    Integration = complete simulation (200 time steps),
    Newton = solving nonlinear system \eqref{eq:system-str-aero},
    Eval UVLM = evaluating $\Maer$ and $\faer$,
    Eval structure = evaluating $\Mstr$ and $\ff\tsp{str}$.}
  \label{tab:wea:solvers}
  \begin{tabular}{l*8r}
    \toprule
    & \mc2c{Exact Newton} & \mc3c{Quasi-Newton} & \mc3c{Inexact Newton} \\
    & \time{00} & \%E & \time0 & \%E & \%Q & \time{00} & \%E & \%I \\
    \midrule
    \mc9c{\it Total} \\[1pt]
    Integration       & 1383.36 &  100
                      &  433.46 & 31.3 &
                      & 1307.81 & 94.5 \\[2pt]
    \N Eval UVLM      &  345.69 &
                      &  346.23 & &
                      &  348.37 & &  \\
    \N Eval structure & 0.14 &
                      & 0.14 & &
                      & 0.14 & & \\
    \N Newton         & 1037.53 &  100
                      &   87.09 &  8.4 & \full
                      &  959.30 & 92.5 & \full \\
    \M Eval UVLM      &  954.32 & 92.0
                      &   76.24 &  7.3 & 87.5
                      &  945.70 & 91.1 & 98.6 \\
    \M Eval structure &    1.85 &  0.2
                      &    1.78 &  0.2 &  2.0
                      &    1.83 &  0.2 &  0.2 \\
    \M Linear solver  &   80.92 &  7.8
                      &    8.93 &  0.9 & 10.3
                      &   11.73 & 11.3 &  1.2 \\[2pt]
    Newton steps      &     992 &&  986 &&& 986 \\
    Refine steps      &     --- &&  --- &&& 986 \\

    \midrule
    \mc9c{\it Average per time step} \\[2pt]
    Integration       &  6.917 &&  2.167 &&&  6.539 \\[2pt]
    \N Eval UVLM      &  1.728 &&  1.731 &&&  1.742 \\
    \N Newton         &  5.188 &&  0.435 &&&  4.797 \\
    \M Eval UVLM      &  4.772 && .38121 &&&  4.729 \\
    \M Eval structure & .00925 && .00891 &&& .00915 \\
    \M Linear solver  & .40462 && .04464 &&& .05865 \\[2pt]
    Newton steps      &  4.935 &&  4.905 &&&  4.905 \\
    \bottomrule
  \end{tabular}
\end{table}

Here the time spent on solving the nonlinear system \eqref{eq:system-str-aero}
dominates the computational effort
in case of the exact and inexact Newton algorithms:
with roughly \SI{1000}{s} (\SI{5}{s} per time step)
it is almost three times as large as the cost for UVLM evaluations
with roughly \SI{350}{s} (\SI{1.7}{s} per time step).
Within the Newton algorithms, however,
the time is still predominantly spent on UVLM evaluations,
which take almost $90$\% or even more for all three algorithms.
The solving of linear systems takes
$7.8$\% (exact), $10.3$\% (quasi-Newton),
and $1.2$\% (inexact), respectively.

As before, the quasi-Newton algorithm is the fastest nonlinear solver
with an even better relative performance:
it is almost 12 times as fast as the exact Newton algorithm,
in part due to the fact that the number of Newton steps is not increased.
In constrast, it is slightly decreased from $992$ to $986$.
This indicates that all approximate inverses are of excellent quality
in this particular simulation.

The inexact Newton algorithm now produces the same iterates
as the quasi-Newton algorithm,
with exactly one linear system solutions per Newton step
(which confirms the excellent quality of all approximate inverses).
Hence the additional cost is immediate from the table:
it amounts to \SI{869.46}{s} of \SI{959.30}{s}
(\SI{4.348}{s} of \SI{4.729}{s} per time step)
for UVLM evaluations
and \SI{2.80}{s} of \SI{11.73}{s}
(\SI{0.014}{s} of \SI{0.0587} per time step)
for checking the residual accuracy in the linear solver.
This is still an extreme overhead between a quasi-Newton algorithm
and the corresponding inexact Newton algorithm
which is based on the quasi-Newton inverse.
In comparison to the exact Newton algorithm, on the other hand,
the cost for UVLM evaluations and linear solving are both reduced,
with the linear solver being responsible
for $89$\% of the total time savings.
This is in line with the typical behavior of inexact Newton techniques.

In summary, the inexact Newton algorithm
now outperforms the exact Newton algorithm by a moderate factor
(and without an artificially increased accuracy)
whereas the quasi-Newton algorithm outperforms both of them
by a factor over ten.


\section{Conclusion}
\label{sec:conclusion}
In this article we have tested
a quasi-Newton type algorithm
and an inexact Newton algorihm
to accelerate nonlinear system solving in implicit integration schemes
for aeroelastic simulations using the unsteady vortex-lattice method.
Both algorithms are based on the structural inverse as approximation
of the full aeroelastic inverse.
On the tested computational examples the structural inverse
approximates the full inverse sufficiently well for
both approximate Newton algorithms to converge.
On these examples the quasi-Newton type algorithm yields a roughly
ten-fold gain in performance over the exact Newton algorithm.
In contrast, the inexact Newton algorithm is found to be
moderately faster than the exact Newton algorithm on the wind turbine problem
but much slower than the quasi-Newton algorithm
although both algorithms generate identical iterates.
On the flexible plate example, the inexact Newton algorithm
is even slower than the exact Newton algorithm
but moderately faster with increased linear termination accuracy.
This strange and unexpected behavior is easily explained by the fact
that the UVLM evaluations dominate the computational cost
to an extent that makes the cost for linear system solving almost irrelevant.
In conclusion, the situation remains unsatisfactory because of two facts.
First, the inexact Newton algorithm is much slower
than the quasi-Newton algorithm
even if both of them generate exactly the same Newton iterates,
just because the latter checks for sufficient accuracy
while the former does not.
Second, the performance of the quasi-Newton algorithm
depends entirely on the approximation quality of the approximate
inverses encountered; the algorithm can diverge
if that quality happens to be very poor,
even in cases where the inexact Newton algorithm
still guarantees quadratic convergence.
Our future work will therefore explore more efficient ways
to compute the derivatives of aerodynamic forces
that are required to check sufficient accuracy of the linear residuals.


\section*{Acknowledgement}

This work would not have been possible without Cristian Guillermo Gebhardt
who initiated the conceptual work and modeling and even the coding
of our aeroelastic simulation approach
and who was a constant driving force of its development.
Further, the authors gratefully acknowledge the financial support from the
Deutsche Forschungsgemeinschaft (DFG, German Research Foundation)
-- SFB1463 -- 434502799.

\bibliographystyle{siam}
\bibliography{literature}

\end{document}